\documentstyle[amssymb,amsfonts,12pt]{amsart}

\newtheorem{theorem}{Theorem}[section]
\newtheorem{lemma}[theorem]{Lemma}
\newtheorem{corollary}[theorem]{Corollary}
\newtheorem{proposition}[theorem]{Proposition}
\newtheorem{definition}[theorem]{Definition}

\newtheorem{conjecture}[theorem]{Conjecture}

\theoremstyle{remark}
\newtheorem{remark}[theorem]{Remark}

\def\QSet{\mbox{\rm\kern.24em
\vrule width.03em height1.48ex depth-.051ex \kern-.26em Q}}

\def\T{{\bf T}}
\def\S{{\bf S}}
\def\Z{{\mathbb Z}}
\def\R{{\mathbb R}}
\def\N{{\mathbb N}}
\def\C{{\mathbb C}}
\def\E{{\mathbb E}}

\def\size{{\operatorname{size}}}

\def\H{{\mathcal H}}

\def\BMO{{\operatorname{BMO}}}
\def\F{{\mathcal F}}

\def\I{{\bf I}}
\def\A{{\mathcal S}}

\def\BMO{{\operatorname{BMO}}}

\def\bas{\begin{align*}}
\def\eas{\end{align*}}
\def\bi{\begin{itemize}}
\def\ei{\end{itemize}}
\newenvironment{proof}{\noindent {\bf Proof} }{\endprf\par}
\def \endprf{\hfill  {\vrule height6pt width6pt depth0pt}\medskip}
\def\emph#1{{\it #1}}

\begin{document}
\title{On some maximal multipliers in $L^p$}

\author{Ciprian Demeter}
\address{Department of Mathematics, Indiana University, 831 East 3rd St., Bloomington IN 47405}
\email{demeterc@@indiana.edu}

\begin{abstract}
We extend an $L^2$ maximal multiplier result of Bourgain to all $L^p$ spaces, $1<p<\infty$.
\end{abstract}

\maketitle

\section{Introduction}

Consider a finite set  $\Lambda=\{\lambda_1,\ldots,\lambda_N\}\subset\R$ and let $D>0$ be the separation constant
$$D:=\min_{1\le i\not=j\le N}|\lambda_i-\lambda_j|.$$ For each $k\in \Z$ with $2^{k}<10^{-2}D$ define $R_k$ to be the collection\footnote{We will abuse notation and also denote by $R_k$ the union of the intervals in the collection} of all intervals of length $2^{k}$ centered at some element from $\Lambda$. The following result was proved in \cite{Bo1}.

\begin{theorem}
\label{gsat2378623890jklcdnvgef}
For each $f\in L^2(\R)$,
\begin{equation}
\label{eq:Bourgain}
\|\sup_{k}|\int_{R_k}\widehat{f}(\xi)e^{2\pi i\xi x}d\xi|\|_{L^2_x}\lesssim (\log N)^2\|f\|_{2}
\end{equation}
\end{theorem}

In \cite{BKO} it was proved that at least a factor of $(\log N)^{1/4}$ is needed on the right hand side. Inequality \eqref{eq:Bourgain} found multiple applications in the work of Bourgain on various ergodic averages, see for example \cite{Bo1}, \cite{Bo15}, \cite{Bo14}. Lacey \cite{La} has used the inequality to prove bounds for the bilinear maximal function. More recently, it has become apparent that variants of \eqref{eq:Bourgain} play a significant role in the analysis of maximal truncations associated with modulation invariant operators, see \cite{D}, \cite{DTT} and \cite{RTT}.

For each $1\le r<\infty$ and each sequence $(x_k)_{k\in \Z}$ in a Hilbert space $\H$, define the \emph{$r$-variational norm} of  $(x_k)_{k \in \Z}$  to be
$$\|x_k\|_{V^r_k(\H)}:=\sup_{k}\|x_k\|_{\H}+ \|x_k\|_{\tilde V^r_k(\H)}$$
where $\tilde V^r_k(\H)$ is the homogeneous $r$-variational seminorm
$$ \|x_k\|_{\tilde V^r_k(\H)}:= \sup_{M,\;k_0<k_1< \ldots <k_{M}}(\sum_{m=1}^M\|x_{k_m}-x_{k_{m-1}}\|_{\H}^r)^{1/r}.$$
If $\H$ is the set of complex numbers $\C$, then the dependence on $\H$ will be suppressed. The only other instance that will occur is when $\H=l^2(N)$ (i.e. $\R^{N}$ with the euclidean inner product). In that case we will use the notation $V^r_k(N)$ and $\tilde V^r_k(N)$. Also, $\|x_Q\|_{V^r_{v_1\le Q\le v_2}}$ will refer to the variational norm of a sequence $x_Q$ indexed by some quantity $Q$ that takes dyadic values between $v_1$ and $v_2$.

For each interval $\omega\in R_k$, let $m_{\omega}$ be a complex valued Schwartz function adapted to $\omega$, that is, supported on $\omega$ and satisfying
$$\|\partial^{\alpha}m_{\omega}\|_{\infty}\le |\omega|^{-\alpha},\;\;\alpha\in\{0,1\}.$$

Let $(w_{\omega})_{\omega\in R_k}$ be a collection of complex weights. Define
$$
\Delta_kf(x):=\sum_{\omega\in R_k}\int w_{\omega}m_{\omega}(\xi)\widehat{f}(\xi)e^{2\pi i\xi x}d\xi,
$$
and also
$$\|w_{\omega}\|_{V^{r,*}}:=\max_{1\le n\le N}\|\{w_{\omega_k}:\lambda_n\in \omega_k\in R_k\}\|_{V^{r}_k},$$
$$\|m_{\omega}\|_{V^{r,*}}:=\max_{1\le n\le N}\|\{m_{\omega_k}(\lambda_n):\lambda_n\in \omega_k\in R_k\}\|_{V^{r}_k}.$$

The following weighted version of \eqref{eq:Bourgain} was proved in \cite{RTT}

\begin{theorem}
\label{BlemmaL2}
For each $r>2$ and $f\in L^2(\R)$ we have the inequality
$$\|\sup_{2^{k}<10^{-2}D}|\Delta_kf(x)|\|_{L^2_x(\R)} \lesssim N^{1/2-1/r}\|w_{\omega}\|_{V^{r,*}}(1+\|m_{\omega}\|_{V^{r,*}})\|f\|_2,$$
with the implicit constant depending only on $r$.
\end{theorem}

Theorem \ref{BlemmaL2} was one of the tools needed to prove the following extension of Bourgain's Return Times theorem:

\begin{theorem}[\cite{RTT}]
\label{Bretthmour}
Let ${\bf X}=(X,\Sigma,\mu, \tau)$ be a dynamical system, and let $1<q\le \infty$ and $p\ge 2$.
For each function $g\in L^{q}(X)$  there is a universal set $X_0\subseteq X$ with $\mu(X_0)=1$, such that for each second dynamical system  ${\bf Y}=(Y,\F,\nu,\sigma)$, each $f\in L^{p}(Y)$  and each $x\in X_0$, the averages
$$\frac1{N}\sum_{n=0}^{N}g(\tau^nx)f(\sigma^ny)$$
converge for $\nu$-almost every $y$.
\end{theorem}

In \cite{RTT}, Theorem \ref{Bretthmour} was proved in the case $p=2$, and the case $p>2$ followed immediately due to the nestedness of $L^p(Y)$ spaces. If one wants to extend the approach from \cite{RTT} to generalize Theorem \ref{Bretthmour} even further, to the case $p<2$, then an $L^p$ version of Theorem \ref{BlemmaL2} needs to be proved. We achieve this here, see Theorem \ref{Blemma} bellow. This result gives hope that the following conjecture can be proved, and it is a first step in the direction of its resolution:

\begin{conjecture}
\label{conjconh}
Theorem \ref{Bretthmour} holds whenever $1<p,q\le \infty$ and
$$\frac1p+\frac1q<\frac32.$$
\end{conjecture}

\begin{theorem}
\label{Blemma}
For each $1<p<2$, each $\epsilon>0$ and each $r>2$ we have the inequality
$$\|\sup_{2^{k}<10^{-2}D}|\Delta_kf(x)|\|_{L^p_x(\R)} \lesssim N^{1/p-1/r+\epsilon}\|w_{\omega}\|_{V^{r,*}}(1+\|m_{\omega}\|_{V^{r,*}})\|f\|_p,$$
with the implicit constant depending only on $r$, $\epsilon$ and $p$.
\end{theorem}

This material is based upon work supported by the National Science Foundation under agreement No. DMS-0635607. In addition, the author was supported by NSF Grant DMS-0556389.  Any opinions, findings and conclusions or recommendations expressed in this material are those of the author and do not necessarily reflect the views of the National Science Foundation.

\section{Where the difficulty lies}

Theorem \ref{Blemma} would typically be applied to the Conjecture \ref{conjconh} with $r$ very close to (but larger than\footnote{The weights $w_\omega$ are typically averages at scale $|\omega|^{-1}$ of the weight function $g$; the $r$ variation of these averages is only bounded when $r>2$}) 2, in which case the exponent of $N$ approaches $1/p-1/2$. On the other hand, this is also the best possible exponent one can hope for in Theorem \ref{Blemma}. This was already observed in \cite{RTT}, and it suffices to consider a single scale. We reproduce the construction for completeness.
\begin{proposition}
For each $N\in\N$ and $p\in (1,2)$ there is a choice of signs $(\varepsilon_n)_{1\le n\le N}$ such that if $\widehat{f}_N=1_{[0,N]}$ then
$$\|\int \widehat{f}_N(\xi)\sum_{l=0}^{N-1}\varepsilon_n1_{[n,n+1]}(\xi)e^{2\pi i\xi x}d\xi\|_{L^p_x(\R)}\gtrsim N^{1/p-1/2}\|f_N\|_{L^p(\R)}.$$
 \end{proposition}
\begin{proof}
 It  immediately follows that  $$\|f_N\|_{L^p(\R)}\sim N^{1-1/p},$$
$$\|(\sum_{n=0}^{N-1}|\int \widehat{f}_N(\xi)1_{[n,n+1]}(\xi)e^{2\pi i\xi x}d\xi|^2)^{1/2}\|_{L^p_x(\R)}\sim N^{1/2}.$$ Khintchine's inequality ends the proof.
\end{proof}

It is worth mentioning that Theorem \ref{Blemma} also holds in the case $p>2$, with the bound\footnote{Note that the argument above also shows that this bound is essentially best possible, when $r$ is close to 2 and $\epsilon$ is close to 0} $N^{\frac12-\frac{2}{pr}+\epsilon}$
but the proof in this case is immediate. Indeed, due to interpolation with the $L^2$ result in Theorem \ref{BlemmaL2}, it suffices to achieve the crude bound $N^{1/2}$ for each $p>2$. To get this, note that by the inequality of Cauchy-Schwartz we have
$$\sup_{2^k\le 10^{-2}D}|\Delta_kf(x)|\le N^{1/2}\|w_{\omega}\|_{V^{r,*}}(\sum_{n=1}^{N}(M(f_n)(x))^2)^{1/2},$$
where
$$f_n(x)=\int_{|\xi-\lambda_n|<D/10} \widehat{f}(\xi)e^{2\pi i\xi}d\xi,$$
and $M(f)$ denotes the Hardy-Littlewood maximal function of $f$. The result then follows from a combination of the Fefferman-Stein inequality \cite{FS} and  Rubio de Francia's result \cite{RdF} on the $L^p$ boundedness ($p>2$) of the square function
$$SQ(f)(x):=(\sum_{n=1}^{N}|f_n(x)|^2)^{1/2}.$$

The main difficulty in proving Theorem \ref{Blemma} is that the square function $SQ(f)$, while it has $L^p\to L^p$ norm independent of $N$ in the case $p\ge 2$, it will have a norm of magnitude $N^{1/p-1/2}$ in the case $p<2$. To overcome this deficit, what we will do instead is relate $\|\sup_{k}|\Delta_kf(x)|\|_{L^p_x(\R)}$ to a more complicated square function, $Vf$, without the use of the Cauchy-Schwartz inequality (and thus without any further loss in powers of $N$). Most of the paper is then devoted to proving that $V$ has essentially the same $L^p,1<p<2$ operator norm as $SQ$ (see Theorem \ref{mainvardghc7836f743fhj}). This itself is a result of independent interest, and  we use time-frequency techniques to prove it. It would be interesting to know whether one could use a more standard approach, via interpolation with an $L^1\to L^{1,\infty}$ result. In particular, it is not clear whether the operator norm $\|V\|_{L^1\to L^{1,\infty}}$ is also of order $N^{1/p-1/2}$. Our approach seems to shed no light on this issue.

We mention that both Theorem \ref{gsat2378623890jklcdnvgef} and Theorem  \ref{BlemmaL2} have variants in \cite{Bo1} and \cite{RTT} respectively, in which the  separation restriction on $\lambda_n$ is eliminated. In each case, the boundedness of the unrestricted operator $\sup_{k\in\Z}|\Delta_kf|$ follows from the boundedness of the restricted version and by using techniques related to the Rademacher-Menshov theorem. A similar approach can be pursued with Theorem \ref{Blemma} here, too. One would first have to bound locally the $L^p$ norm by the $L^2$ norm and then apply Rademacher-Menshov type arguments to the localized operator. These local contributions are then summed over the whole real line, to produce the desired estimate. We refer the interested reader to the proof of Theorem 8.7 in \cite{RTT} for details.

\section{A variational inequality}
The following simple inequalities will be useful throughout the rest of the paper.

\begin{lemma}
\label{lem:auxq3}
For each $k$, let $a_k,b_k$ be some complex numbers. Then for each $r\ge 1$
$$\|a_kb_k\|_{V^r_k}\lesssim \|a_{k}\|_{V^r_k}\|b_{k}\|_{V^r_k},$$
\end{lemma}

\begin{lemma}
\label{lem:auxq3kkt}
Let $\{c_k:=(c_k^{(n)})_{1\le n\le N}:\;k\in\Z\}\subseteq l^2(N)$. Then
$$\|c_k\|_{V^r_k(N)}\le (\sum_{n=1}^{N}\|c_k^{(n)}\|_{V^r_k}^2)^{1/2}.$$
\end{lemma}

\begin{lemma}
\label{lem:auxq3kkthgdhgd}
Let $(x_k)_{k_0\le k< k_L}\in\C$ and let $k_0<k_1<\ldots k_L$. Then
$$\|x_k\|_{V^r_{k_0\le k< k_L}}\le \sum_{l=0}^{L-1}\|x_k\|_{V^r_{k_l\le k< k_{l+1}}}.$$
\end{lemma}

There is an extensive literature on variational estimates in Harmonic Analysis and Ergodic Theory. We confine ourselves to mentioning only a few such papers: \cite{CJKRW},\cite{JKRW},\cite{JSW}.

For each\footnote{We can ignore the dyadic points, since they have measure zero} $x\in\R$ and $k\in \Z$ denote by $I(x,k)$  the unique dyadic interval of length $2^{-k}$ which contains $x$. Let $\phi^{(k)}$ be Schwartz functions with  Fourier transform adapted to $[-\frac12,\frac12]$. For each interval $I$ of length $2^{-k}$, define
$$\phi_I(y):=\frac{1}{|I|}\phi^{(k)}(\frac{y-c(I)}{|I|}).$$

\begin{lemma}
\label{varlemma}
We have for each $1<p<\infty$ and each $r>2$
$$\left\|\|\langle f, \phi_{I(k,x)}\rangle\|_{V^r_k}\right\|_{L^p_x}\lesssim (1+\|\widehat{\phi^{(k)}}(0)\|_{V^r_{k}})\|f\|_{p},$$
with the implicit constant depending only on $\phi$, $r$ and $p$.
\end{lemma}
\begin{proof}
Consider the dyadic martingale
$$\E_k(f)(x)=\langle f, \frac{1}{|I(k,x)|}1_{I(k,x)}\rangle.$$
A result of Lepingle \cite{Lep} implies that
$$\left\|\|\E_k(f)(x)\|_{V^r_k}\right\|_{L^p_x}\lesssim \|f\|_{p}.$$
By using this and Lemma \ref{lem:auxq3} we get that
$$\left\|\|\widehat{\phi^{(k)}}(0)\E_k(f)(x)\|_{V^r_k}\right\|_{L^p_x}\lesssim \|\widehat{\phi^{(k)}}(0)\|_{V^r_{k}}\|f\|_{p}.$$

It suffices now to note that the difference operator can be bounded as follows
\begin{align*}
\left\|\|\langle f, \frac{\widehat{\phi^{(k)}}(0)}{|I(k,x)|}1_{I(k,x)}-\phi_{I(k,x)}\rangle\|_{V^r_k}\right\|_{L^p_x}&\le \left\|\left(\sum_{k}|\langle f, \frac{\widehat{\phi^{(k)}}(0)}{|I(k,x)|}1_{I(k,x)}-\phi_{I(k,x)}\rangle|^2\right)^{1/2}\right\|_{L^p_x}\\&\lesssim \|f\|_{p}.
\end{align*}
The first inequality holds since $r>2$, while the second one is a consequence of the Littlewood-Paley theory, since $\phi_I-\frac{\widehat{\phi^{(k)}}(0)}{|I|}1_I$ has mean zero (when $|I|=2^{-k}$).
\end{proof}

\section{Time-frequency interlude}
The purpose of this section is to prepare the ground for the proof of Theorem \ref{mainvardghc7836f743fhj} in the next section.

To prove Theorem \ref{Blemma}, it will suffice by scaling invariance to assume that
\begin{equation}
\label{sepcond}
\min_{1\le i\not= j\le N}|\lambda_i-\lambda_j|= 1.
\end{equation}
We will do so throughout the rest of the paper.

\begin{definition}
A tile $s=s(I,n)$ will be a rectangle of the form $$s(I,n):=I\times [\lambda_n-2^{k-1},\lambda_n+2^{k-1}],$$
where $1\le n\le N$, while $I$ is a dyadic interval of length $2^{-k}>100$. We refer to $I_s:=I$ and $\omega_s:=[\lambda_n-2^{k-1},\lambda_n+2^{k-1}]$ as the time interval and frequency interval of the tile $s$.
\end{definition}

The collection of all tiles will be denoted by $\S$. We define $\S_{n}$ to consist of the collection of all tiles $s\in \S$ such that $\lambda_n\in \omega_s$. Note that $(\S_{n})_{n=1}^{N}$ forms a partition of $\S$.

\begin{definition}
A collection $\S'\subset \S$ will be referred to as convex, if whenever $s,s''\in\S'$, $s'\in\S$, $\lambda_n\in \omega_s\cap\omega_{s'}\cap\omega_{s''}$ for some $n$ and  $I_s\subseteq I_{s'}\subseteq I_{s''}$, these also imply that $s'\in\S'$.
\end{definition}

\begin{definition}
A tree $(\T,T)$ with top $T\in \S$ is a collection of tiles $\T\subset \S$ such that
$I_s\subseteq I_T$ and $\omega_T\subseteq \omega_s$ for each $s\in \T$.
\end{definition}

Note that each tree is entirely contained in some $\S_{n}$.

We will denote by $T_m$ the Fourier projection associated with the multiplier $m$:
$$T_mf(x):=\int \widehat{f}(\xi)m(\xi)e^{2\pi i\xi x}d\xi.$$

We will use the notation
$$\tilde{\chi}_{I}(x)=(1+\frac{|x-c(I)|}{|I|})^{-1}.$$

\begin{definition}
\label{adsfdfdwqryetluidhweilh}
Let $f$ be a $L^2$ function  and let $(\T,T)$ be a tree.
We define the size $\size(\T)$ of $\T$ relative\footnote{The function with respect to which the size is computed will change throughout the paper; however, it will always be clear from the context} to $f$ as
$$\size(\T):=\sup_{s\in \T}\sup_{m_s}\frac{1}{|I_s|^{1/2}}\|\tilde{\chi}_{I_s}^{10}(x)T_{m_s}f(x)\|_{L^2_{x}},$$
where  $m_s$ ranges over all functions adapted to $10\omega_s$.
\end{definition}

For a fixed tile $s$, the quantity
$$\sup_{m_s}\frac{1}{|I_s|^{1/2}}\|\tilde{\chi}_{I_s}^{10}(x)T_{m_s}f(x)\|_{L^2_{x}},$$
is an approximate measure for the $L^2$ norm of the portion of $f$ that is localized in time-frequency in $s$.

We need some notation and results from \cite{MTT8}.

\begin{definition}\label{hull-definition}
Let $(\T,T)$ be a tree, and let $\I_T$ be the collection of all
maximal dyadic intervals $I\subseteq I_T$ which have the property that
$3 I$ does not contain any of the intervals $I_s$ with $s\in T$.
For an integer $j$ with $2^{-j}\le |I_T|$ let $\tilde{E}_j$ be the union of all intervals $I$
in $\I_T$ such that $|I|< 2^{-j}$

For an integer $j$ with $2^{-j}>|I_T|$ we define $\tilde{E}_j=\emptyset$.
\end{definition}
The sets $\tilde{E}_j$ obviously depend on the tree $\T$, but we suppress
this dependence.

We recall the following Lemma 4.12 from \cite{MTT8}.
\begin{lemma}
\label{count-again}

Let $\Omega_j$ be the collection of connected components of $\tilde E_j$. Then $\Omega_j$ is a finite collection of dyadic intervals each of which has length equal to an integer multiple of $2^{-j}$.

For each $I \in \Omega_j$, let $x^l_I$ and $x^r_I$ denote the left and right endpoints of $I$, and let $I^l_j$ and $I^r_j$ denote the intervals
$$ I^l_j := (x^l_I - 2^{-j-1}, x^l_I - 2^{-j - 2})$$
$$ I^r_j := (x^r_I + 2^{-j - 2}, x^r_I + 2^{-j-1}).$$
Then the intervals $I^l_j$ are disjoint as $j$ varies in the integers
with $2^{-j}\le |I_T|$ and $I$ varies in $\Omega_j$.

Similar statements hold for the $I^r_j$.
\end{lemma}
Recall also the weight function from \cite{MTT8}.
$$\mu_{j,\T}(x):=\sum_{j}2^{-|j'-j|/100}\sum_{y\in \partial \tilde {E}_j}(1+2^{j'}|x-y|)^{-100}.$$

An important role in our future investigation is played by the function
$$W_\T(x):=\sum_{s\in\T}\frac{1_{I_s}(x)}{|I_s|}\int \mu_{j(s)}(y)\tilde{\chi}_{I_s}^2(y)dy,$$
where $|I_s|=2^{-j(s)}$.

In the same Lemma 4.12 from \cite{MTT8} it was proved (and this is an immediate consequence of Lemma \ref{count-again}) that

$$\sum_{j}2^{-j}\# \partial \tilde {E}_j \lesssim |I_T|.$$
This in turn automatically implies that
$$\int W_\T(x)dx\lesssim |I_T|.$$
We will need a stronger result, namely that
$$\int_J W_\T(x)dx\lesssim |J|$$
for each dyadic $J\subseteq I_T$. This however also turns out to be an easy consequence of Lemma \ref{count-again}. We leave the details to the reader. The following is then a standard corollary.

\begin{corollary}
\label{corerrrrror}
We have
$$\|W_\T\|_{\BMO}\lesssim 1,$$
where $\|\cdot\|_{\BMO}$ stands for the dyadic BMO norm.
\end{corollary}

The following result is a simplified version\footnote{One of the simplifications arises from the simpler definition of size from Definition \ref{adsfdfdwqryetluidhweilh}} of Proposition 7.4 from \cite{MTT8}. In its current form, it
appears essentially in \cite{DeTh}.

Each tree $\T$ defines a region in the phase-space domain. The operator $V$ from Theorem \ref{mainvardghc7836f743fhj} will receive a contribution from all the trees. The lemma below will show that for each tree, the most important contribution comes from the part of the function $f$, denoted by $\Pi_\T(f)$, that essentially lives in the region defined by $\T$.

\begin{lemma}[Phase-space projections]
\label{phase-spaceproj34}
Let $(\T,T)$ be a tree and $f\in L^2(\R)$\footnote{Since we only want an $L^2$ estimate in (ii), there is no need for the extra assumption $\|f\|_{\infty}\le 1$ from \cite{DeTh} or \cite{MTT8}}. There is a function $\Pi_{\T}(f)$, called the phase-space projection of $f$ on the tree $\T$, satisfying the properties:
\begin{itemize}
\item (i) (Control by size) $\Pi_{\T}(f)$ is supported in $2I_T$ and satisfies the $L^{\infty}$ bound
$$\|\Pi_{\T}(f)\|_{\infty}\lesssim \size(\T).$$
As a consequence, for each $1\le t\le \infty$,
\begin{equation}
\label{contrbysize}
\|\Pi_{\T}(f)\|_{t}\lesssim |I_T|^{1/t}\size(\T)
\end{equation}
with the implicit constant independent of  $\T$ or $f$.

\item (ii) ($\Pi_{\T}(f)$ approximates $f$ on $\T$) For each $s\in\T$ with $|\omega_s|=2^{j}$, each $m$ adapted to $\omega_s$ we have
\begin{equation}
\label{contrbysize1}
\|(|I_s|\phi_{I_s})^{1/2}T_m(\Pi_{\T}(f)-f)\|_{2}\lesssim \size(\T)|I_s|^{-1/2}\int \tilde{\chi}_{I_s}^2(x) \mu_j(x)dx.
\end{equation}

\end{itemize}
\end{lemma}

Assume that for each $k\in \Z$ we have  $n$ functions $\phi_1^{(k)},...,\phi_n^{(k)}$ whose Fourier transforms are adapted to the interval $[-1/2,1/2]$.
Define also the modulated versions
\begin{equation}
\label{modulatedphi}
\phi_{I,n}(x):=\frac{1}{|I|}\phi_n^{(k)}(\frac{x-c(I)}{|I|})e^{2\pi i\lambda_n x},
\end{equation}
whenever $|I|=2^{-k}$.
Also, if $s\in\S_n$, define $\phi_s:=\phi_{I_s,n}$. We will assume throughout the rest of the section that
$$\max_{n}\|\widehat{\phi_{n}^{(k)}}(0)\|_{V^r_{k,2^{-k}>100}}\le 1.$$

If $(\T,T)$ is a tree and $r>2$, define
$$V_\T f(x):=\|1_{I_s}(x)\langle f, \phi_{s}\rangle \|_{V^r_{s\in\T}}.$$

As a consequence of Lemma \ref{varlemma}, Corollary \ref{corerrrrror} and Lemma \ref{phase-spaceproj34} we have
\begin{proposition}
\label{dbcdsghf3opkqfdnvjerk}
Let $(\T,T)$ be a tree and let $f\in L^2(\R)\cap L^t(\R)$, for some $1<t<\infty$. Then
\begin{equation}
\label{error6734767}
\|V_\T f\|_t\lesssim |I_T|^{1/t}\size(\T).
\end{equation}
\end{proposition}

\begin{proof}
Triangle's inequality implies
\begin{equation*}
\|V_\T f\|_t \le \|V_\T\Pi_\T(f)\|_t+\|V_\T(f-\Pi_\T(f))\|_t.
\end{equation*}

Since the first term above is well controlled by Lemma \ref{varlemma} and by \eqref{contrbysize}, it suffices to control the second term. Note that by  \eqref{contrbysize1} and the fact that $\widehat{\phi}_s$ is supported in $\omega_s$, we have for each $s\in\T$

\begin{align*}
|\langle f-\Pi_\T(f), \phi_s\rangle|&\le \int|\phi_{I_s}(x)T_{m_s}(f-\Pi_\T(f))(x)|dx\\&\lesssim \|(|I_s|\phi_{I_s}(x))^{1/2}T_{m_s}(f-\Pi_\T(f))(x)\|_2\|(|I_s|^{-1}\phi_{I_s})^{1/2}\|_2\\&\lesssim \size(\T)|I_s|^{-1}\int \tilde{\chi}_{I_s}^2(x) \mu_{j(s)}(x)dx,
\end{align*}
where $m_s$ is an appropriate function, adapted to $\omega_s$ and equal to 1 on the support of $\widehat{\phi}_s$.

Next, by using the above and the trivial estimate $\|x_k\|_{V^r_k}\le \|x_k\|_{l^1(\Z)}$, we obtain the following pointwise estimate

$$|V_\T(f-\Pi_\T(f))(x)|\le \size(\T)W_\T(x).$$
The claim is now immediate from Corollary \ref{corerrrrror} and John-Nirenberg's inequality.
\end{proof}

The following lemma is immediate.

\begin{lemma}
\label{Linfestforsize777}
If $(\T,T)$ is a tree and let $f\in L^2(\R)\cap L^{\infty}(\R)$ then
$$\size(\T)\lesssim \sup_{s\in\T}\inf_{x\in I_s}Mf(x).$$
\end{lemma}

We now show how to split $\S$ into trees with special properties.

\begin{lemma}
\label{Bessel777}
Let $\S'\subset \S$ be a finite collection of tiles and let $f\in L^2(\R)\cap L^{\infty}(\R)$. Then we can split
$$\S'=\bigcup_{2^{-m}\le 2\size(\S')}\A_m,$$
where $\A_m$ is the union of disjoint trees
$$\A_m=\bigcup_{\T\in\F_m}\T$$
with $\size(\T)\le 2^{-m}$ for each $\T\in\F_m$ and moreover
\begin{equation}
\label{bcsdhgc34or843;fre0-956]-3-u9}
\sum_{\T\in\F_m}|I_\T|\lesssim 2^{2m}\|f\|_2^2.
\end{equation}
\end{lemma}

This lemma is obtained by iterating the following lemma

\begin{lemma}
\label{Bessel777hhhh}
Let $\S'\subset \S$ be a finite collection of tiles, and assume
$$\size(\S')\le 2\lambda.$$
Then one can split
$$\S'=\S'_1\cup\S'_2,$$
where $\S'_1$ is the union of disjoint trees
$$\S'_1=\bigcup_{\T\in\F}\T$$
with
$$\sum_{\T\in\F}|I_\T|\lesssim \lambda^{-2}\|f\|_2^2,$$
while
$$\size(\S'_{2})\le \lambda.$$
\end{lemma}

\begin{proof}
The proof is entirely standard. We briefly sketch the argument.
\begin{itemize}
\item Step 0: Set $\A:=\S'$ and $\F:=\emptyset$
\item Step 1: Select a tile $T\in \A$ with maximal time interval $I_T$ such that there exists $m_T$ adapted to $10\omega_T$ satisfying

\begin{equation}
\label{hgdweghuifyuewhcjkdsbc}
\frac{1}{|I_T|^{1/2}}\|\tilde{\chi}_{I_T}^{10}(x)T_{m_T}f(x)\|_{L^2_{x}}>\lambda.\end{equation}

If there is more than one  tile having a fixed (maximal) time interval that
qualifies to be selected, just select any of them. Construct the maximal tree $\T$ in $\A$ having top $T$. Reset $\F:=\F\cup\{\T\}$ and  $\A:=\A\setminus \T$.

If there is no such tile to be selected, stop, the algorithm is over.

\item Step 2: Go to Step 1
\end{itemize}
After the algorithm stops, set $\S'_1:=\bigcup_{\T\in\F}\T$ and $\S'_2:=\S'\setminus \S'_1$, and note that $\size(\S'_2)\le \lambda$, as desired.

Observe also that the rectangles $(I_T\times 10\omega_T)_{\T\in \F}$ are pairwise disjoint. Indeed, if $\T,\T'$ have the tops $T$ and $T'$ in the same $\S_n$ ($1\le n\le N$) then $I_T$ and $I_{T'}$ will be disjoint by maximality. If $T$ and $T'$ belong to distinct $\S_n$ and $\S_{n'}$, then, while their time intervals may intersect, the frequency components $10\omega_T, 10\omega_{T'}$ will be disjoint, due to \eqref{sepcond} and the fact that $|\omega_T|,|\omega_{T'}|<\frac{1}{100}$.

Finally, inequality \eqref{hgdweghuifyuewhcjkdsbc} combined with this pairwise disjointness implies via a (now) very standard $TT^{*}$ argument (see for example Corollary 7.6 in \cite{MTT}) that

$$\sum_{\T\in\F}|I_T|\lesssim \lambda^{-2}\|f\|_2^2.$$

\end{proof}

\begin{remark}
Note that if $\S'$ is convex, then each $\A_m$ is convex.
\end{remark}

\section{A square function estimate for the variational norm}

Let as before $\phi_1^{(k)},...,\phi_n^{(k)}$ be Schwartz functions whose Fourier transforms are adapted to the interval $[-1/2,1/2]$, and let $\phi_{I,n}$ be defined as in \eqref{modulatedphi}.

For $r>2$ define
$$Vf(x):=(\sum_{n=1}^N\|1_I(x)\langle f,\phi_{I,n}\rangle\|_{V^r_{I,|I|>100}}^2)^{1/2}.$$

\begin{theorem}
\label{mainvardghc7836f743fhj}
We have for each $1<p\le 2$ and each $\epsilon>0$
\begin{equation}
\label{hgdchgerfuo7834890rukvjhfjvh}
\|Vf\|_{p}\lesssim N^{\frac1p-\frac12+\epsilon}(1+\max_{n}\|\widehat{\phi_{n}^{(k)}}(0)\|_{V^r_{k,2^{-k}>100}})\|f\|_p,
\end{equation}
with the implicit constant depending only on $r$, $p$ and $\epsilon$.
\end{theorem}
\begin{proof}
We can and will assume that $\max_{n}\|\widehat{\phi_{n}^{(k)}}(0)\|_{V^r_{k,2^{-k}>100}}\le 1$. The inequality when $p=2$ (with $\epsilon=0$) follows immediately from Lemma \ref{varlemma} and by orthogonality. By interpolating with the $L^2$ bound, it suffices to prove \eqref{hgdchgerfuo7834890rukvjhfjvh} with the bound $N^{1/2}$ rather than $N^{\frac1p-\frac12+\epsilon}$

We rephrase the inequality to be proved in terms of tiles. We are to prove that
$$V_{\S'}f(x):=(\sum_{n=1}^N\|1_{I_s}(x)\langle f,\phi_{s}\rangle\|_{V^r_{s\in\S'_n}}^2)^{1/2}$$
satisfies
\begin{equation*}
\|V_{\S'}f\|_{p}\lesssim N^{1/2}\|f\|_p,
\end{equation*}
uniformly over all finite $\S'\subseteq \S$ which in addition are also convex.

Also, by restricted type interpolation, it suffices to prove
\begin{equation}
\label{hgdchgerfuo7834890rukvjhfjvhjkvhjkfhvj}
|\{x\in\R: V_{\S'}1_{F}(x)>\lambda^{1-\epsilon} N^{1/2}\}|\lesssim \frac{|F|}{\lambda^p},
\end{equation}
for each $F\subset \R$ with finite measure and each $\epsilon>0$, with an implicit constant depending just on $p$ and $\epsilon$. Moreover, due to the $L^2$ result, it suffices to prove \eqref{hgdchgerfuo7834890rukvjhfjvhjkvhjkfhvj} for $\lambda<1$.

Consider the exceptional set
$$E:=\{x\in\R: M1_F(x)\ge C\lambda\},$$
with $C$ large enough, independent of $F$.
Since $|E|\le 2C^{-p}\frac{|F|}{\lambda^p},$ it further suffices to prove
$$|\{x\in\R\setminus E: V_{\S'}1_{F}(x)>\lambda^{1-\epsilon} N^{1/2}\}|\lesssim \frac{|F|}{\lambda^p}.$$

Obviously, on $\R\setminus E$ only the tiles $s\in\S'$ with $I_s\cap (\R\setminus E)\not=\emptyset$ will contribute to $V_{\S'}1_{F}$. Thus we can assume $\S'$ consists only of such tiles. Note that $\S'$ remains convex. Lemma \ref{Linfestforsize777} implies now that $\size(\S')\le \lambda/2$, if $C$ is sufficiently large. Here, the size is computed with respect to the function $1_F$.

Recall now the decomposition in Lemma \ref{Bessel777}. For each $2^{-m}\le \lambda$ and each $\T\in\F_m$ we define the exceptional set
$$E_\T:=\{x\in I_T: V_{\T}1_F(x)>\lambda^{ \frac12-\epsilon}2^{-m/2}\},$$
and we note that due to Proposition \ref{dbcdsghf3opkqfdnvjerk},
$$|E_\T|\lesssim \lambda^{-t(\frac12-\epsilon)}2^{-mt/2}|I_\T|.$$
Furthermore, due to \eqref{bcsdhgc34or843;fre0-956]-3-u9} we get that
$$|E':=\bigcup_{2^{-m}\le \lambda}\bigcup_{\T\in \F_m}E_\T|\lesssim \sum_{2^{-m}\le \lambda}\lambda^{-t(\frac12-\epsilon)}2^{-mt/2}2^{2m}|F|\lesssim \frac{|F|}{\lambda^p},$$
assuming $t>\max\{4,\frac{2-p}{\epsilon}\}$.
\end{proof}
As a consequence, it further suffices to prove that for $x\in\R\setminus (E\cup E')$ we have
$$|V_{\S'}1_{F}(x)|\le \lambda^{1-\epsilon} N^{1/2}.$$

Observe that each $x\in\R$ belongs to at most  $N$ time intervals $I_T$ of tops $T$ of trees in $\A_m$, for each $2^{-m}\le \lambda$. For each $m$, at most one of these trees is in each $\S_n$ ($1\le n\le N$). Call these trees $\T_{1,m},\ldots,\T_{n,m}$, where we allow some of them to be empty. Then
$$V_{\A_m}1_{F}(x)=(\sum_{n=1}^{N}(V_{\T_{n,m}}1_{F}(x))^2)^{1/2}.$$

By the triangle inequality in Lemma \ref{lem:auxq3kkthgdhgd} (possible since each $\A_m$ is convex) combined with the triangle inequality for the norm of $l^2(N)$, it follows that whenever $x\in\R\setminus (E\cup E')$ we have
\begin{align*}
|V_{\S'}1_{F}(x)|&\le \sum_{2^{-m}\le \lambda}|V_{\A_m}1_{F}(x)|\\&\le \sum_{2^{-m}\le \lambda}N^{1/2}\lambda^{ \frac12-\epsilon}2^{-m/2}\\&=\lambda^{1-\epsilon} N^{1/2}.
\end{align*}

\section{Proof of Theorem \ref{Blemma}}
\label{s:ghhghg}

We first use Fourier series to rewrite the maximal function $\sup_k|\Delta_kf|$ as  a maximal expression involving sums of exponentials with variable coefficients. We then estimate  the $L^p$ norm  by the $L^2$ norm, on each dyadic interval of unit length. The $L^2$ norm of these maximal exponential sums will then be  estimated by using the metric entropy method, see Lemma \ref{lem:auxq2} below. The result is that the local contribution, $\|\sup_k|\Delta_kf|\|_{L^p(J)}$, will be bounded by  $\|Vf\|_{L^p(J)}$. These local contributions are then summed up and  Theorem \ref{mainvardghc7836f743fhj} finishes the argument.

The idea of reducing matters to maximal exponential sums goes back to \cite{Bo1}. The advantage of dealing with expressions of the form
$$\sup_{k}|\sum_{n=1}^{N}c_{k,n}e^{2\pi i \lambda_n y}|$$
resides both in the almost orthogonality of the exponentials, but also in the finite dimensionality of the coefficient space, $l^2(N)$. In \cite{Bo1}, the reduction to exponential sums is performed via an averaging procedure, that relies crucially on Plancherel's theorem. However, in our case, this type of argument is not available, since we deal with $L^p$ estimates. We will instead use windowed Fourier series expansions. The key is the fact that we can reduce matters to $L^2$ locally, without further loss in the constants.

The following lemma  essentially appears in \cite{Bo1}. Its current formulation is taken from \cite{RTT}.

\begin{lemma}[Lemma 8.4,\;\cite{RTT}]
\label{lem:auxq2}
For each  set $C=\{c_k\}\subseteq l^2(N)$ and each $r>2$ we have
$$\|\sup_{k}|\sum_{n=1}^{N}c_{k,n}e^{2\pi i \lambda_n y}|\|_{L^2_y([0,D^{-1}))}\lesssim N^{1/2-1/r}\|c_k\|_{V^r_k(N)},$$
with the implicit constant depending only on $r$.
\end{lemma}

\begin{proof}[of Theorem \ref{Blemma}]
By rescaling, we can assume that the separation constant $D=1$.

Let $\psi:\R\to \C$ be a smooth function such that $\widehat{\psi}$ is supported on $[-1,1]$ and equals 1 on [-1/2,1/2]. For each dyadic interval $I:=[l/2^k,(l+1)/2^k]$ denote by $$\phi_{I,n}(x):=\frac{1}{|I|}\widehat m_{\omega,*}(\frac{x-c(I)}{|I|})e^{2\pi i\lambda_n x},$$
where $\omega$ is the unique interval in $R_k$ that contains $\lambda_n$, $m_{\omega,*}$ is the rescaled, shifted at the origin version of $m_{\omega}$, that is
$$m_{\omega,*}(\frac{\xi-\lambda_n}{|\omega|})=m_{\omega}(\xi).$$

Assume now $\omega$ is centered at some $\lambda_n$. By using windowed Fourier series on $\omega$ we can write
$$
\widehat f(\xi)m_{\omega}(\xi)=\frac{1}{|\omega|}\sum_{l\in \Z}[\int \widehat f(\eta)m_{\omega}(\eta)e^{2\pi i \frac{l}{|\omega|}\eta}d\eta]e^{-2\pi i \frac{l}{|\omega|}\xi}{\widehat \psi}(\frac{\xi-\lambda_n}{|\omega|}),$$
so that
$$\int m_{\omega}(\xi)\widehat{f}(\xi)e^{2\pi i\xi x}d\xi=\sum_{|I|=|\omega|^{-1}\atop{I \;dyadic}}e^{2\pi i\lambda_n x}\langle f, \phi_{I,n}\rangle \psi_I(x).
$$
where $$\psi_I(x):=\psi(\frac{x-c(I)}{|I|}).$$

Let now $J$ be a dyadic interval of unit length. For each $2^k<\frac{1}{100}$ and each $l\in \Z$ denote by $I(J,k,l)$ the unique dyadic interval of length $2^{-k}$ such that $l2^{-k}+I(J,k,l)$ contains $J$. We have
\begin{align*}
\|\sup_{k}|\Delta_kf(x)|\|_{L^p_x(J)}&\le \sum_{l\in\Z}\|\sup_{k}|(\sum_{n\le N}e^{2\pi i\lambda_n x}\langle f, \phi_{I(J,k,l),n}\rangle w_{k,n})\psi_{I(J,k,l)}(x)|\|_{L^p_x(J)}\\&\lesssim \sum_{l\in\Z}2^{-100|l|}\|\sup_{k}|\sum_{n\le N}e^{2\pi i\lambda_n x}\langle f, \phi_{I(J,k,l),n}\rangle w_{k,n}|\|_{L^p_x(J)},
\end{align*}
where $w_{k,n}:=w_{\omega}$, given that $\omega$ is centered at $\lambda_n$ and has scale $2^{k}$. By invoking standard rescaling arguments it suffices to analyze the case $l=0$, so any further dependence on $l$ will be suppressed. We can now invoke Lemma \ref{lem:auxq2},  Lemma \ref{lem:auxq3kkt} and Lemma \ref{lem:auxq3} to get
\begin{align*}
\|\sup_{k}|\sum_{n\le N}e^{2\pi i\lambda_n x}\langle f, \phi_{I(J,k),n}\rangle w_{k,n}|\|_{L^p_x(J)}&\le \|\sup_{k}|\sum_{n\le N}e^{2\pi i\lambda_n x}\langle f, \phi_{I(J,k),n}\rangle w_{k,n}|\|_{L^2_x(J)}\\&\lesssim N^{1/2-1/r}\max _{n}\|w_{k,n}\|_{V_k^r}(\sum_{n=1}^N\|\langle f, \phi_{I(J,k),n}\rangle\|_{V_k^r}^2 )^{1/2}\\&=N^{1/2-1/r}\|w_\omega\|_{V^{r,*}}V_Jf,
\end{align*}
where $V_Jf$ denotes the (constant) value of $Vf(x)$ on $J$.

By raising to the $p^{th}$ power and by summing up over all dyadic $J$ of unit length, it suffices now to prove that
\begin{equation}
\|Vf\|_{L^p(\R)}\lesssim N^{\frac1p-\frac12+\epsilon}(1+\|m_{\omega}\|_{V^{r,*}})\|f\|_{L^p(\R)}.
\end{equation}
This however was proved in Theorem \ref{mainvardghc7836f743fhj}.
\end{proof}


\end{document}